\documentclass{article}
\usepackage[margin=25mm]{geometry}
\usepackage[utf8]{inputenc}
\usepackage{authblk}
\usepackage[shortlabels]{enumitem}

\usepackage{hyperref}
\usepackage{enumitem} 
\newlist{condenum}{enumerate}{1} 
\setlist[condenum]{label=\bfseries P\arabic*., 
                   ref=\arabic*, wide}

\usepackage{algpseudocode}
\usepackage{%
	amsfonts,%
	amssymb,%
	amsthm,%
	bbm,%
	float,%
	etoolbox,%
	mathrsfs,%
	mathtools,%
	multirow,%
	pgf,%
	subfigure,%
	tikz,%
}
\usepackage[normalem]{ulem}
\usepackage{cleveref}
\pgfdeclareplotmark{mystar}{
    \node[star,star point ratio=2.25,minimum size=6pt,
          inner sep=0pt,draw=black,solid,fill=red] {};
}
\pgfdeclareplotmark{mystar2}{
    \node[star,star point ratio=2.25,minimum size=6pt,
          inner sep=0pt,draw=black,solid,fill=blue] {};
}

\newlist{temp}{enumerate}{1} 
\setlist[temp]{label=(\roman*), 
                   ref=\roman*, wide}
\def\ind{\mathbbm{1}}
\usepackage{bbm}
\usepackage{breqn}

\newcommand\numeq[1]%
  {\stackrel{\scriptscriptstyle(\mkern-1.5mu#1\mkern-1.5mu)}{=}}

\theoremstyle{plain}
\newtheorem{thm}{Theorem}
\newtheorem{lem}[thm]{Lemma}

\newtheorem{cor}[thm]{Corollary}
\newtheorem{remark}{Remark}

\theoremstyle{definition}

\newtheoremstyle{case}{}{}{}{}{}{:}{ }{}
\theoremstyle{case}

\theoremstyle{remark}

\DeclarePairedDelimiter\abs{\lvert}{\rvert}%

\newcommand{\brac}[1]{\left(#1\right)}
\newcommand{\cbrac}[1]{\left\{#1\right\}}
\newcommand{\sbrac}[1]{\left[#1\right]}

\newcommand{\expect}[1]{\mathbb{E}\sbrac{#1}}

\providecommand{\keywords}[1]
{
  \small	
  \textbf{\textit{Keywords---}} #1
}

\title{On the Performance of Large Loss Systems with Adaptive Multiserver Jobs}

\author[1]{Samira Ghanbarian}
\author[2]{Arpan Mukhopadhyay}
\author[3]{Fabrice M. Guillemin}
\author[4]{Ravi R. Mazumdar}
\affil[1,4]{Department of Electrical and Computer Engineering, University of Waterloo}
\affil[2]{Department of Computer Science, University of Warwick}
\affil[3]{Orange Innovation}

\date{} 

\begin{document}
\maketitle

\begin{abstract}
In this paper, we study systems where each job or request can be split into a flexible number of sub-jobs up to a maximum limit. 
The number of sub-jobs a job is split into depends on the number of available servers found upon its arrival.
All sub-jobs of a job are then processed in parallel at different servers leading to a linear speed-up of the job. 
We refer to such jobs as {\em adaptive multi-server jobs}. A vast majority of jobs submitted to modern data centres can be modelled as adaptive multi-server jobs. These jobs can also represent requests for files in a file-server system where each file is stored at multiple locations from where different parts of the file can be downloaded in parallel. We study the problem of optimal assignment of such requests  when each server can process at most one sub-job at any given instant and there is no waiting room in the system.

We assume that, upon arrival, a job can only access a randomly sampled subset of $k(n)$ servers from a total of $n$ servers, and the number of sub-jobs is determined based on the number of idle servers within the sampled subset. We analyze the steady-state performance of the system when system load varies according to $\lambda(n)  =1 - \beta n^{-\alpha}$ for $\alpha \in [0,1)$, and $\beta \geq 0$. Our interest is to find how large the subset $k(n)$ should be in order to have zero blocking and maximum speed-up in the limit as $n \to \infty$. We first characterize the system's performance when the jobs have access to the full system, i.e., $k(n)=n$. In this setting, we show that the blocking probability approaches to zero at the rate $O(1/\sqrt{n})$ and the mean response time of accepted jobs approaches to its minimum achievable value at rate $O(1/n)$. We then consider the case where the jobs only have access to subset of servers, i.e., $k(n) < n$. We show that  as long as $k(n)=\omega(n^\alpha)$, the same asymptotic performance can be achieved as in the case with full system access. In particular, for $k(n)=\Theta(n^\alpha \log n)$, we show that both the blocking probability and the mean response time approach to their desired limits at rate $O(n^{-(1-\alpha)/2})$. We use a combination of Stein's and Lyapunov drift methods to establish state space collapse and characterize the performance of the system for large $n$.

\end{abstract}


\maketitle

\section{Introduction}

Consider a system consisting of $n$ parallel servers, each working at a unit rate. Jobs arrive at the system according to a Poisson process and bring exponentially distributed workloads. There is no queueing buffer in the system and jobs that cannot find an available server immediately upon arrival, are blocked. Jobs are parallelizable. That is, each job can split into smaller parts, named sub-jobs or tasks, where each sub-job is processed in parallel with the other sub-jobs resulting in a liner speed-up in its execution. Thus, a job can run on multiple servers at the same time. Compared with traditional one-server jobs where each job can only join a single server, multiserver parallelizable jobs benefit from the simultaneous use of multiple resources and hence have a better execution time. The question we address in this paper is under what conditions can we achieve the maximum gain in the system performance, in terms of average response time and blocking probability? We address this question for situations when jobs have access to the set of full servers or a limited subset of sampled servers. Ideally, we would like the number of sampled servers to be small compared to the total number of servers available in the system.

Classical one-server jobs do not represent the complex structure of modern requests anymore. A large number of jobs in modern data centers are multiserver parallelizable jobs, which means they can run on more than one server in parallel. For instance, many $5G/6G$ network functions are decomposed into a large number of smaller microservices in parallel in order to improve the network latency \cite{Fabrice_microservice}. Google's Borg system \cite{Google_Borg_2015}, schedules each job to a set of parallel processors for better performance and high utilization. MapReduce framework \cite{MapReduce_2008} and machine learning jobs such as TensorFlow \cite{TensorFlow_2016}, allow parallelized computation across large clusters of servers to reduce the network traffic. Coded storage systems such as Erasure Codes \cite{ErasureCodes_2017}, are another example where an original file is recovered from multiple data chunks stored in parallel.

Most multiserver parallelizable jobs are moldable \cite{moldable_2003, moldable_2002}. This means that these jobs are flexible in the number of servers they can join and they can adapt based on the availability of resources. This is in contrast to multiserver jobs that require a fixed and predetermined number of servers \cite{Brill_multiserver_1984, FILIPPOPOULOs_multiserver_2007, Zero_wait_central-2021, sharp_zer_central_2022}. In the latter case, if the number of available servers is insufficient for the job's requirements, the job will not start processing, resulting in a waste of resources. Therefore, adaptive multiserver jobs allow for better utilization of system resources, but at the cost of increased system complexity.

While there is some flexibility in the number of servers at which a job can be processed in parallel and obtain a linear speed-up, achieving full flexibility is normally not possible. Most computing jobs only benefit from parallelization as long as the number of servers used for parallel processing is below a certain limit; beyond this limit it may not be possible to further break down a job or breaking it further will only result in sub-linear speed-up resulting from the increased overhead of parallelization~\cite{Berg_speedup_2020}.

Therefore, in this work, we focus on adaptive multiserver jobs which can be processed by at most $d$ servers in parallel. If $1 \leq l\leq d$ servers are used to process a job, then a linear speed-up of $l$ is obtained. If no empty server is found upon the entry of a job, then the job is discarded or {\em blocked}.
Our goal in this paper is to design job assignment schemes which result in the maximum average speed-up of $d$ per job while maintaining zero blocking in the steady-state. A natural job assignment scheme to use in this scenario is the {\em greedy} scheme where all empty servers, up to a maximum of $d$ servers, is used to process each incoming job. The first question we address in this paper is whether this greedy scheme is able to achieve the desired objective of zero blocking and maximum speed-up as the system size becomes large ($n \to \infty$).

In a real system, it may not always be possible for a job to have access to the states of all the servers. 
For example, consider a system with $m$ dispatchers and $n$ servers where jobs arrive at a rate that scales as $O(n)$. A token based scheme in which the servers send their idleness status as tokens to each of the $m$ dispatchers would require  a message rate of $O(nm)$ messages per unit time and $O(n)$ memory size at each server if the scheme were to perform optimally~\cite{gamarnik2018delay}. This messaging and memory overhead can quickly become prohibitively high when the number $m$ of dispatchers is $O(n)$ as is often the case in modern large-scale data centres~\cite{lu2011join}. 
To account for this fact, we consider a model where, upon arrival, a job only has access to a random subset of servers of size $k(n)$, where $k(n) \ll n$. This is equivalent to the situation where each arriving job is randomly assigned to one of the $m$ dispatchers and the dispatcher then sampling $k(n)$ servers at random from the system. We consider the same greedy scheme as discussed before, i.e., the dispatcher assigns the job to a maximum of $d$ servers found available among the sampled subset of $k(n)$ servers. Note that this scheme is easy to implement even with a large number of dispatchers as the messaging rate is only $O(n k(n))$ messages per unit time for the same arrival rate and the memory required at each dispatcher is only $O(k(n))$.
Therefore, our second objective in the paper is to determine 
if it is possible to achieve the desired objective of zero blocking and maximum speed-up with only small value of $k(n)$, i.e., by knowing the states of only a small number of servers at each arrival instant.



\subsection{Contributions} \label{section: contributions}

The model we consider is a generalization of the multi-rate loss model that has been studied in the telecommunications literature \cite{Kaufman, Whitt_dropping_1985}. The novelty of this model lies in two main aspects: Firstly, the number of servers used to process a job is not preset and depends on the state of the sampled system upon arrival. Secondly, the service times are now dependent on the number of servers assigned to a job. This is in contrast to multi-rate loss models, where each class requires a pre-determined number of servers and the service time distribution is the same for all jobs of a given class and not dependent on the state of the system. This generalization introduces sojourn time as a performance measure that depends on the number of parallel servers used to process the job. In other words, the model we consider is a loss model with state-dependent service times where the dependence is based on the state seen by an arrival. An interesting facet of the model is that there is no simple coupling or dominance possible due to these properties since the service time durations change depending on the number of available servers at the arrival time. Our main contributions are summarized as follows:

\begin{itemize}
    \item We show that when the greedy scheme is used and each job has access to all the servers, i.e., when $k(n) = n$, the number of jobs that are processed over fewer than $d$ servers eventually vanish in the stationary regime as the system size grows. Furthermore, all arriving jobs tend to find at least $d$ available servers in the system. We establish the above properties for a job arrival rate of $n\lambda(n)$, where $\lambda(n)$ varies as $\lambda(n)=1-\beta n^{-\alpha}$ for $\alpha \in [0,1)$ which includes the so called {\em mean-field regime} ($\alpha=0$), the {sub-Halfin-Whitt} regime ($\alpha \in (0, 0.5)$), the {Halfin-Whitt} regime ($\alpha=0.5$), and the {\em super-Halfin-Whitt regime} ($\alpha \in (0.5,1)$). 
    Specifically, our results show that the blocking probability approaches to zero and the mean response time approaches to its minimum possible value of $1/d$ at rates $O(1/\sqrt{n})$ and $O(1/n)$, respectively, in all the regimes mentioned above.


    \item  In the case where each job has access to only a subset of size $k(n) \ll n$ sampled uniformly at random, we show that the greedy scheme still achieves the same performance asymptotically as long as the size $k(n)$ of the sampled subset of servers grows as $k(n)=\omega(n^\alpha)$. In particular, when $k(n)=\Theta(n^\alpha \log n)$, we show that both the mean response time and the blocking probability exhibit a convergence to their respective optimal values at the rate $O(n^{-\brac{1-\alpha}/2})$ where $\alpha$ denotes the rate at which the arrival rate of jobs approaches to the critical load of the system. This implies that when the load varies as $\lambda(n)=1-\beta n^{-\alpha}$, it is possible to achieve asymptotically optimal performance with a message rate of only  $O(n^{1+\alpha} \log n)$ per unit time and a memory size of $O(n^{\alpha} \log n)$ per dispatcher. This is a significant improvement over the $O(n^2)$ message rate and $O(n)$ memory size per dispatcher required for join-the-idle-queue-type schemes with $O(n)$ dispatchers.

    \item In order to prove our results, we use a combination of Stein's method, Lyapunov drift method, and sample-path arguments~\cite{liu2020steady,braverman2020steady}. To elaborate more, in our system, it is difficult to characterize the exact fluid limit due to the dependence of $k(n)$ and $\lambda(n)$ on the system size. Instead, we use a simple deterministic dynamical system which mimics the desired system behaviour for large values of $n$. By comparing the generator of the original system and to that of the deterministic system, we obtain bounds on the quantities of interest. The Lyapunov functions that we use in our analysis are carefully constructed to take into account the fact that the service times of jobs now depend on the number of available servers found upon entry.
    
\end{itemize}


\subsection{Related Work} \label{section: related work}

The study of parallelizable jobs has been of great interest to researchers. Different queueing systems are introduced to model parallelizable jobs. One classic theoretical model is Fork-Join (FJ) system. In an FJ system with $n$ servers, an incoming job splits into $n$ independent tasks at the fork station, and each task joins one server, processed in a First Come First Serve order. When all tasks are complete, the delay of the job is determined by the maximum delay of its tasks at the join station. The stationary joint workload distribution of FJ systems with two servers is studied in \cite{FJ_two_servers_1984}. However, the analysis of FJ systems in the general case of more than two servers is extremely challenging due to the complex correlations between the fork and join stations. The existing literature only provides approximations and bounds on the performance metrics of the system such as upper bounds on the delay performance of the system, as can be seen in studies like \cite{baccelli_FJ_1989} and \cite{Rizk_FJ_2016}. A comprehensive survey of the existing results on FJ systems is presented in \cite{FJ_survey}.

A generalization of FJ systems is partial or limited FJ systems, where in a system with $n$ servers, the incoming job splits into $ k<n$ tasks. In the case of fixed $k$, upper bounds on the tail distribution of the response time of jobs are studied in \cite{Rizk_FJ_2016}, when tasks are assigned to servers randomly. The study is extended in \cite{heterogeneous_partial_FJ_2023} to heterogeneous systems with slow and fast servers, where probabilistic selection policies for task assignment are considered, achieving upper bounds on average completion time. Additionally, for a varying number of tasks $k(n)$, it is shown in \cite{limited_FJ_Wang_as_indep_2019} that any $k(n) =o(n^{1/4})$ subset of servers becomes asymptotically independent that leads to an upper bound on the mean response time of jobs. Furthermore, in \cite{Wang_distributed_2020}, a batch-filling policy for parallelizable jobs is proposed, where tasks are assigned sequentially to the shortest queue in a sampled set of servers. This policy achieves zero queueing delay, meaning that every task of the job starts its process immediately.

In limited FJ systems, each server has its own queue, and tasks wait in the queue of the server to start their service. However, an alternative queueing model is when there is a central queue from which tasks are dispatched to servers. These systems are called multiserver-job systems. If the number of available servers is less than the number of job's tasks, the job will wait in the central queue until there are enough resources available. The stationary distribution of the number of jobs in systems with two servers is studied in \cite{Brill_multiserver_1984} and \cite{FILIPPOPOULOs_multiserver_2007}. In multiserver-job systems, the job at the head of the queue may not fit into the system immediately and therefore incurs an extra waiting time. However, it is shown in \cite{Zero_wait_central-2021} and \cite{sharp_zer_central_2022} as the system load, the number of tasks, and the number of servers scale, jobs can achieve zero asymptotic waiting time. Analogous results are studied in \cite{multiserver_asymptotic_opimality_2022}, where the mean response time is minimized under heavy traffic limit, while the number of servers remains unscaled. An alternative approach to address the extra waiting time is by dropping the job if it cannot fit into the system. Multiserver-job systems with blocking have been of considerable interest and are studied well such as in \cite{Whitt_dropping_1985}, where the goal is to minimize the blocking probability of each job.

Another class of queueing models for parallelizable jobs involves variable server allocation. Prior research in this domain has predominantly focused on the use of speedup functions; these functions measure the ratio of a job's original service time to its accelerated service time when parallelized. In other words, speedup functions denote the amount of acceleration jobs can get by joining multiple servers. One well-known example of a speedup function is Amdahl's law \cite{Amdahl's_law_2008}, which accounts for scenarios where only a part of the job is parallelizable, while the remaining part receives no parallelization. In \cite{Berg_speedup_2018}, the optimality of concave and sublinear speedup functions for exponential job sizes is studied. It is shown that the policy that shares servers equally among the jobs in the system minimizes the mean response time. In \cite{Berg_speedup_2020}, it is assumed that some jobs can split into any number of parts, while others follow a threshold parallelization approach where the level of parallelization is limited to a specific threshold. Optimal policies are derived for this speedup function when job sizes are exponential. Additionally, in \cite{WCFS_Grososf_2022}, a work-conserving finite skip framework is introduced, which includes threshold parallelism as a special case. The mean response time of the system is characterized under heavy traffic conditions.

The study most closely related to our model is the threshold speedup function, as discussed in \cite{WCFS_Grososf_2022}. However, there are several significant distinctions between our model and the aforementioned study. Firstly, they consider a queueing model while our model functions as a blocking system, where jobs are blocked if they cannot find available resources. Secondly, in \cite{WCFS_Grososf_2022}, it is assumed that all servers are accessible to all jobs and there is no sampling upon arrival. Conversely, we consider both cases of sampling and non-sampling upon arrival.

\subsection{General Notation and Organization}

The following notation will be used throughout the paper. Let $f(n)$ and $g(n)$ be positive real-valued increasing functions. We define $f(n)$ to be $o\brac{g(n)}$ if $\sup \lim_{n \to \infty}\frac{f(n)}{g(n)} = 0$, $f(n)$ to be $O\brac{g(n)}$ if $\sup \lim_{n \to \infty}\frac{f(n)}{g(n)} < \infty$, $f(n)$ to be $\omega\brac{g(n)}$ if $\inf \lim_{n \to \infty}\frac{f(n)}{g(n)} = \infty$, and $f(n)$ to be $\Omega\brac{g(n)}$ if $\inf \lim_{n \to \infty}\frac{f(n)}{g(n)}>0$. Moreover, we define $f(n)$ to be $\Theta\brac{g(n)}$ if $f(n)$ is $O\brac{g(n)}$ and $\Omega\brac{g(n)}$. Additionally, for an integer $m$, we use $[m]$ to represent the set $\{1,2,\ldots,m\}$.

The rest of the paper is organized as follows. In Section~\ref{section: system model}, we introduce the system model. Main results on the asymptotic optimality of the system and its convergence rates for $k(n)=n$ are presented in Section~\ref{section: full server access}. These results are extended to limited system access with $k(n)<n$ in Section~\ref{section: Limited server access}. Finally, concluding remarks and potential extension directions are discussed in Section~\ref{section: conclusion}.

\section{System Model} 
\label{section: system model}
We consider a system with $n$ servers. Each server can process only one job at any given time and there is no waiting room in the system. It is assumed that each job can be processed simultaneously at a maximum of $d\geq 1$ servers. A job's inherent processing time, i.e., the time it would take to process the job at one server, is assumed to be independent and exponentially distributed with unit mean. When a job is run in parallel on $i \in [d]$ servers, a speed-up of $i$ is obtained, i.e, the job's processing time is reduced by a factor of $i$ in comparison to its inherent service time.

\begin{remark}
When a job joins $i \in [d]$ idle servers, it is partitioned into $i$ equal smaller components termed sub-jobs or tasks. The completion of a job occurs when all of its constituent tasks have left the system.
Since all sub-jobs have equal sizes and join the idle servers simultaneously, their processing time is identical. This synchronized execution ensures that all sub-jobs are complete at the same time, resulting in the processing time for the entire job being equivalent to that of a single sub-job. Consequently, for a job with size $\mathcal{E}$ and subject to a maximum degree of parallelization $d$, the minimum attainable response time is $\frac{\mathcal{E}}{d}$.
\end{remark}

Jobs arrive at the system according to a Poisson process with the rate $n\lambda (n)$, where $\lambda (n) = 1- \beta n ^{-\alpha}\geq 0$ for some $\alpha \in [0,1)$ and $\beta>0$. Varying the value of $\alpha$ enables the study of the system under different traffic regimes:
(i) $\alpha=0, \beta \in (0,1)$ corresponds to the {\em mean-field regime}, (ii) $\alpha=1/2$ corresponds to the {\em Halfin-Whitt regime}, and (iii) $\alpha \in (0,1/2)$ (resp. $\alpha \in (1/2,1)$)  corresponds to the {\em sub-Halfin-Whitt (resp. super-Halfin-Whitt) regime}.

When a job arrives, it must be immediately and irrevocably assigned to a subset of at most $d$ servers where it is processed until it leaves the system. In this paper, our main focus is to study the performance of the system under the following two scenarios: 

\begin{enumerate}
    \item {\bf Full system access}: In this scenario, each job has access to the states of all the servers in the system.

    \item {\bf Limited system access}: In this scenario, each job has access to the states of only a subset of $k(n)$ servers chosen uniformly at random from $n$ servers.
\end{enumerate}
Note that $k(n)=n$ corresponds to the case with full system access. A job must be assigned to at most $d$ of the $k(n)$ servers the job has access to. If none of these $k(n)$ servers are found idle upon the entry of the job, then the job is discarded or {\em blocked}, which corresponds to a loss. Our goal is to design job assignment schemes which achieve minimum mean response time of jobs while maintaining zero blocking in the steady-state.

To define the problem more concretely, we introduce the following notations.
For each $i \in [d]$, we let $X_i(t)$ denote the number of jobs that are being processed  simultaneously at $i$ servers  at time $t\geq 0$. 
Define $x_i(t), i\in [d],$ as $x_i(t)=X_i(t)/n$, i.e, the scaled number of jobs being processed simultaneously at $i$ servers at time time $t$. Clearly, under any Markovian job assignment scheme (a scheme which makes assignments based on the current state of the system), the process $x(\cdot)=(x_i(\cdot), i\in [d])$ is Markov with a unique stationary distribution. By omitting the explicit dependence on $t$, we denote by $x=(x_i, i\in [d])$ the state of the system distributed according to the stationary distribution. Additionally, we define the fraction of busy servers at stedy-state as $q_1 = q_1(x) = \sum_{i \in [d]} ix_i$, and the fraction of idle servers at steady-state as $q_0 = 1 - q_1$.

Let $D$ denote the time spent by an accepted job in the system
in the steady-state and $P_b$ denote the steady-state blocking probability of a job. Since the total departure rate of jobs at steady-state is $n\expect{\sum_{i \in [d]}i x_i} = n \expect{q_1}$, applying the rate conservation principle and Little's law, we obtain

\begin{align}
    \lambda(n) (1-P_b)&=\expect{q_1},\label{eq:rate_conservation}\\
    \lambda(n)(1-P_b)\expect{D}&=\expect{\sum_{i\in [d]}x_i}\label{eq:little}.
\end{align}

Since a job can be processed at a maximum of $d$ servers simultaneously, the minimum possible value of the average response time of jobs, $\expect{D}$, is $1/d$. From the above equations, it is clear that for any Markovian job assignment scheme to achieve the minimum steady-state mean response time of $\frac{1}{d}$ while maintaining a zero blocking at steady-state (i.e., $P_b=0$), $\expect{x_i}$ for each $i\in [d]$ must satisfy:

\begin{align}
\expect{x_i}&=0, \quad \forall i\in [d-1],\\
\expect{x_d}&=\frac{\lambda (n)}{d}.
\end{align}
Let $x^*=\brac{0,\ldots,0,\frac{\lambda(n)}{d}}$. Hence, any scheme for which $\abs{\expect{x_i}-x_i^*} \to 0$ as $n \to \infty$ is asymptotically optimal. Below we define a greedy job assignment scheme which aims to achieve this.




{\bf Greedy Assignment Scheme}: Under the greedy assignment scheme, upon arrival of a job,
if $l\geq 1$ servers are found available among the $k(n)$ servers the job has access to, then $\min(d,l)$ of these available servers are used to process the job, i.e., under the greedy scheme all the available servers up to a maximum of $d$ servers are used to process the incoming job.

In the following subsections, we characterise the performance of the greedy assignment scheme.

\subsection{Comparison with a fluid limit: The Stein's Approach} \label{section: preliminary results}

In order to characterise the system's performance under the greedy scheme, we follow the Stein's method to compare the dynamics of the system under the greedy scheme to that of a deterministic fluid limit. Our aim is to show that in the limit as $n\to \infty$ the system under the greedy scheme essentially behaves as the fluid limit.

{\bf The fluid limit}: The deterministic fluid limit to which we compare the dynamics of the original system is defined through the following system of ODEs:

\begin{align}
\dot x_i&=-ix_i, \quad \forall i\in [d-1],\label{eq:ode1}\\
\dot x_d&=\lambda (n)-dx_d\label{eq:ode2}.
\end{align}
Intuitively, in the above system, all arriving jobs find $d$ or more available servers upon entry. It is easy to see that starting from any initial state the above system converges to $x^*$ as $t \to \infty$. Below we attempt to bound the distance between the steady-state performance of the original system and that of the fluid limit.

To this end, we denote by $A_i(x)$ the probability with which an incoming job is processed at $i \in \{0,1,\ldots,d\}$ servers when the system is in state $x$. Specifically, $A_i(x)$ represents the probability of finding exactly $i$ idle servers when $i \in \cbrac{0,1,\ldots,d-1}$, and the probability of having $d$ or more idle servers when $i=d$. It is important to note that in this context, $A_0(x)$ corresponds to the blocking probability in state $x$, and $\sum_{j=0}^d A_j(x)=1$.

When $k(n)=n$, by the definition of the greedy scheme, these probabilities are given by

\begin{equation}\label{eq: arrival prob case n}
    A_i(x)=\begin{cases}
    &\ind(nq_0=i), \quad \text{ if } i\in \{0,1,2,\ldots,d-1\},\\
    &\ind(nq_0\geq d),\quad  \text{ if } i=d.
    \end{cases}
\end{equation}
Similarly, when $k(n)< n$, the probabilities $A_i(x)$ are given by means of a hypergeometric distribution. Specifically, 

\begin{equation} \label{eq: arrival prob case k}
    A_i(x)=\begin{cases}
    &\frac{\binom{nq_0}{i}\binom{nq_1}{k(n)-i}}{\binom{n}{k(n)}}, \quad \text{ if } i\in \cbrac{0,1,2,\ldots,d-1},\\
    &\sum_{i\geq d}\frac{\binom{nq_0}{i}\binom{nq_1}{k(n)-i}}{\binom{n}{k(n)}}, \quad \text{ if } i=d.
    \end{cases}
\end{equation}


In the lemma below, using Stein's method of generator comparison, we express the mean squared distance between $x_d$ and $x_d^*$ as a function of the probability $A_d(x)$ of finding at least $d$ available servers upon entry at the steady state of the system.
This expression is later used to bound the mean squared distance under different scenarios of server access.

\begin{lem} \label{lemma: general policy lemma}
Under the equilibrium measure of the system:
\begin{equation} \label{eq: first bound for general policy}
\expect{(\lambda (n)-dx_d)^2} =  \frac{d}{n} \lambda(n)\expect{A_d(x)}+\lambda (n)\expect{(1-A_d(x))(\lambda (n)-dx_d)}.
\end{equation}
In particular, when $k(n)=n$, we have
\begin{equation} \label{eq: first bound for general policy, k=n}
\expect{(\lambda (n)-dx_d)^2} =  \frac{d}{n} \lambda(n) \expect{A_d(x)}+\expect{\ind\brac{q_1> 1-\frac{d}{n}}(\lambda (n)-dx_d)}.
\end{equation}
\end{lem}

\begin{proof}
We consider the Lyapunov function $V(x)=\frac{1}{2d}(\lambda (n)-dx_d)^2$. Let $G$ be the generator of the Markov chain $x(\cdot)$. Furthermore, let $L$ be the generator of the system of ODEs given by~\eqref{eq:ode1}-\eqref{eq:ode2}.
The crux of the Stein's approach is to compare the drift of the function $V$ under $G$ to that under $L$. We note that under $L$ the drift of $V$ is given by
\begin{equation}
LV(x) = \frac{\partial V}{\partial x_d}(x)\dot{x}_d= -(\lambda (n)-dx_d)\dot x_d= -(\lambda (n)-dx_d)^2.
\end{equation}
Similarly, under $G$ the drift of $V$ is given by

\begin{equation*}
    GV(x) = \sum_{x^' \not = x} r(x,x^') \brac{V(x^')-V(x)},
\end{equation*}
where $r(x,x^')$ denotes the transition rate from state $x$ to state $x^'$. Given that the probability of having $d$ or more idle servers is $A_d(x)$, the function $V(x)$ will transit from state $x$ to state $x + \frac{e_d}{n}$ at the rate $n\lambda(n) A_d(x)$, where $e_d$ denotes the $d$-dimensional unit vector with one at the $d^{\textrm{th}}$ position. Additionally, the system will transit from state $x$ to the state $x-\frac{e_d}{n}$, if one of the jobs occupying $d$ servers departs the system. As there are $n x_d$ jobs split into $d$ sub-jobs, the departures occur at the rate $n d x_d$. Thus,
\begin{equation}
GV(x)=n\lambda (n) A_d(x)\brac{V\brac{x+\frac{1}{n}e_d}-V(x)}+ndx_d\brac{V\brac{x-\frac{1}{n}e_d}-V(x)}.
\end{equation}

Under the equilibrium measure,  $\expect{GV(x)}=0$, and hence
\begin{align}
\expect{GV(x)-LV(x)}=\expect{-LV(x)}=\expect{(\lambda (n)-dx_d)^2}.
\end{align}
Therefore, to bound $\expect{(\lambda (n)-dx_d)^2}$ it is sufficient to bound $\expect{GV(x)-LV(x)}$.
This is done as follows. Using Taylor series expansion of $V$ and noting that $L V(x)=-\brac{\frac{\partial V}{\partial x_d}}^2$ we have
\begin{align}
\expect{GV(x)-LV(x)}
&=\expect{n\lambda (n) A_d(x)\brac{\frac{1}{n}\frac{\partial V}{\partial x_d}(x)+\frac{1}{2n^2}\frac{\partial^2 V}{\partial x_d^2}(\xi)}}\nonumber\\
&+\expect{nd x_d\brac{-\frac{1}{n}\frac{\partial V}{\partial x_d}(x)+\frac{1}{2n^2}\frac{\partial^2 V}{\partial x_d^2}(\theta)}}+\expect{\brac{\frac{\partial V}{\partial x_d}(x)}^2},
\end{align}
where $\xi$ and $\theta$ are $d$-dimensional vectors. 
Simplifying the RHS of the above and using $\frac{\partial^2 V}{\partial x_d^2}(y)=d$ for any vector $y$ yields
\begin{align}
\expect{GV(x)-LV(x)}
&=\expect{\brac{\lambda (n) A_d(x)-dx_d+\frac{\partial V}{\partial x_d}(x)}\brac{\frac{\partial V}{\partial x_d}(x)}}+\frac{d}{2n}\expect{\lambda (n) A_d(x)+d x_d}.
\end{align}
Jobs that occupy $d$ servers in the system have an arrival rate of $n\lambda(n) A_d(x)$ and a departure rate of $ndx_d$. Since the system is in a steady state, the rate conservation law applies, requiring these quantities to be equal on average. Hence $\expect{dx_d} = \lambda(n) \expect{A_d(x)}$. Finally, noting that $\frac{\partial V}{\partial x_d}(x)=-(\lambda (n)-dx_d)$ in the above, we have the desired result.
\end{proof}

The first terms on the RHS of Equations~\eqref{eq: first bound for general policy} and \eqref{eq: first bound for general policy, k=n} are $O(\frac{1}{n})$. Therefore, to establish the convergence of $x_d$ to $x_d^*$, it is sufficient to show that the second term in the above expressions  scales as $O(\frac{1}{n})$. In the following sections, we establish this both when the jobs have access to all servers (i.e., $k(n)=n$) and when the jobs have access to only a subset of servers (i.e., $k(n) < n$).

\section{Full System Access} \label{section: full server access}

We begin by considering the case when $k(n)=n$, i.e, when every arrival has access to the complete set of servers. We first show that for $k(n)=n$, the second term in~\eqref{eq: first bound for general policy, k=n}  is always negative.
We use sample path arguments to establish this in the following lemma.

\begin{lem} \label{lemma: k=n second term bound}
For $k(n)=n$, if the system starts at a state where $\sum_{i \in [d-1]}x_i(0)=0$, then at all times we have $\sum_{i\in [d-1]}x_i(t) \leq 1/n$. Furthermore, if $\lambda(n) = 1 - \beta /n ^{\alpha} \geq 0$ for $\alpha \in [0,1)$ and $\beta >0$, then at steady-state we have
\begin{equation}
\expect{\ind\brac{q_1> 1-\frac{d}{n}}(\lambda (n)-dx_d)} \leq 0,
\end{equation}
for all sufficiently large $n$.
\end{lem}

\begin{proof}
Let us denote the sum $\sum_{i \in [d-1]}x_i$
at any state $x$ by $S_{d-1}(x)$.
If the system starts at a state where $S_{d-1}(x)=0$, then $S_{d-1}(x)$ remains zero
until the number of free servers in the system drops strictly below $d$ since all arrivals finding $d$ or more free servers will only increase $x_d$ keeping the other components of $x$ the same. Once the number of free servers drops strictly below $d$,
the next arrival increases  $S_{d-1}(x)$ from zero to at most $1/n$. Let us call this job as the {\em tagged job}. If, upon arrival of the tagged job, the sum $S_{d-1}(x)$ increases to $1/n$, then system becomes fully busy after the arrival of the tagged job as the tagged job occupies all remaining servers. Hence, until the tagged job leaves the system, subsequent arrivals either find the system fully busy (hence are blocked) or find at least $d$ servers available (which occurs if a job occupying $d$ servers departs before the arrival occurs and the tagged job leaves the system). In either case, the $S_{d-1}(x)$ remains constant at $1/n$ until the tagged job departs. 
When the tagged job departs, $S_{d-1}(x)$ again drops to zero. From this point onward we can apply the same chain of arguments as above. This shows that $S_{d-1}(x)$ never increases beyond $1/n$ which proves the first part of the lemma.

If the system starts at a state satisfying $S_{d-1}(x) > 1/n$, then it is easy to see that with non-zero probability the system enters the set of states where $S_{d-1}(x) =0$ in finite time. 
However, the first part of the lemma implies that starting from states satisfying $S_{d-1}(x)=0$, it is not possible to reach states satisfying $S_{d-1}(x) > 1/n$ and the chain always remains in states satisfying $S_{d-1}(x)\leq 1/n$. Hence, the states with $S_{d-1}(x) > 1/n$ are transient and the states with $S_{d-1}(x) \leq 1/n$ form a single, finite communicating class. 
Therefore, there exists a unique invariant probability measure of the chain concentrated only on the states satisfying $S_{d-1}(x) \leq 1/n$.
This implies that at steady-state we have $q_1(x)=\sum_{i \in [d]} i x_i < dx_d+d S_{d-1}(x)\leq dx_d+d/n$ with probability one. Hence, when $q_1 > 1-d/n$, we have $dx_d > 1-2d/n$. When $\lambda (n)$ is a function of $n$ and varies as $\lambda (n)=1-\beta/n^{\alpha}$ for $\alpha \in [0,1)$, it decreases at a faster rate than $1-d/2n$. Thus, it holds that $dx_d > \lambda (n)$ for all $n$ sufficiently large. This shows that 
$\expect{\ind\brac{q_1> 1-\frac{d}{n}}(\lambda (n)-dx_d)} \leq 0$ for all sufficiently large $n$.
\end{proof}

Combining Lemmas~\ref{lemma: general policy lemma} and \ref{lemma: k=n second term bound}, we conclude the following.

\begin{cor} \label{cor: k=n}
Let $k(n)=n$ and $\lambda (n) = 1-\beta/n^{\alpha} \geq 0$ for $\alpha \in [0,1)$ and $\beta >0$. Then, under the equilibrium measure, we have
\begin{equation}\label{eq: k=n, dxd=lambda}
\expect{(\lambda (n)-dx_d)^2}\leq \frac{d\lambda (n)}{n},
\end{equation}
for all sufficiently large $n$. 
\end{cor}

The above Corollary establishes that as the system size $n$ increases, the higher dimensional system undergoes a state space collapse and reduces to a lower dimension. Specifically, the $d$-dimensional system simplifies to one dimension for large system size, where it can be fully described by only the jobs occupying $d$ servers simultaneously. This result has important implications, as it ensures every arrival in the system achieves the minimum possible average response time in the asymptotic limit. Consequently, the system demonstrates asymptotic optimality in terms of average response time.

In the following theorem, we further demonstrate this result, along with the asymptotic optimality in terms of the blocking probability. Specifically, we establish that the system achieves the minimum average response time of $\frac{1}{d}$ at a rate of $O(\frac{1}{n})$ and achieves zero blocking probability at a rate of $O(\frac{1}{\sqrt{n}})$.

\begin{thm} \label{theorem: k=n, performance}
    Let $k(n)=n$ and $\lambda (n) = 1-\beta/n^{\alpha}\geq 0$, where $\alpha \in [0,1)$ and $\beta >0$. Then in the steady-state regime, the blocking probability of the system converges to zero at a rate of $O(\frac{1}{\sqrt{n}})$, and the mean response time of accepted jobs converges to $\frac{1}{d}$ at a rate of $O(\frac{1}{n})$.
\end{thm}    
    \begin{proof}
        By the rate conservation law, we have 
        $$\expect{q_1} = \lambda(n)\brac{1- P_b},$$
        where $P_b$ denotes the blocking probability of the system. On the other hand, from the definition of $q_1 = \sum_{i \in [d]} ix_i$, we have $dx_d \leq q_1$ which leads to $\expect{dx_d} \leq \lambda(n) \brac{1- P_b}$. Hence,

    \begin{equation*}
       P_b^2 \leq \frac{1}{\lambda(n)^2}\brac{\expect{\lambda(n) - dx_d}}^2 \leq \frac{1}{\lambda(n)^2}\expect{\brac{\lambda(n) - dx_d}^2} \leq \frac{d}{n \lambda(n)},
    \end{equation*}
where the second inequality follows from Jensen's inequality and the last inequality follows from corollary~\ref{cor: k=n}. 

Note that $\lambda(n) = 1 - \beta/ n^\alpha$. If $\alpha = 0$ and $0 < \beta < 1$, then it readily follows that 

\begin{equation}\label{eq: k=n, alpha=0, bound on block P}
    P_b^2 \leq\frac{d}{n(1-\beta)},
\end{equation} 
and the blocking probability of the system converges to zero at rate $O(\frac{1}{\sqrt{n}})$. If $\alpha >0$ and $\beta>0$, since $\frac{\beta}{n^{\alpha}}<1$ for large $n$, we have

\begin{equation} \label{eq: k=n, alpha>0, bound on block P}
    P_b^2 \leq \frac{d}{n}\brac{1 + \frac{\beta}{n^{\alpha}}+ O(\frac{1}{n^{2\alpha}})} \leq \frac{d}{n} + o(\frac{1}{n}).
\end{equation}

Hence, again the blocking probability of the system converges to zero at rate $O(\frac{1}{\sqrt{n}})$. 

Consider the mean response time of the system given by $\expect{D}$. From Little's law for the stationary regime, we have

\begin{equation*}
    \lambda(n) (1-P_b) \expect{D}  = \expect{\sum_{i}x_i}.
\end{equation*}

From Lemma~\ref{lemma: k=n second term bound}, the first $d-1$ states of the system form a single finite communicating class in the steady-state, such that $\sum_{i \in [d-1]} x_i \leq 1/n$. Moreover, $dx_d \leq q_1$, due to the definition of $q_1$. Thus

\begin{equation*}
    \expect{D} \leq \frac{1}{\lambda(n) (1-P_b) } \brac{\frac{1}{n} + \frac{\expect{q_1}}{d}}.
\end{equation*}

But we know that $\expect{q_1} = \lambda(n)\brac{1- P_b}$ by the rate conservation law. Therefore,

\begin{equation*}
    \expect{D} \leq \frac{1}{d} + \frac{1}{n \lambda(n) (1-P_b)}.
\end{equation*}

If $\alpha=0$ and $0 < \beta < 1$, then from Equation~\eqref{eq: k=n, alpha=0, bound on block P}, and the fact $\sqrt{\frac{d}{n(1-\beta)}}<1$ for large $n$, we have

\begin{align}
    \expect{D} &\leq \frac{1}{d} + \frac{1}{n(1-\beta)} \brac{1 + \sqrt{\frac{d}{n(1-\beta)}} + O(\frac{1}{n})}
    \nonumber
    \\ & \leq \frac{1}{d} + \frac{1}{n(1-\beta)}+ o(\frac{1}{n}).
\end{align}

If $\alpha>0$ and $\beta >0$, then from Equation~\eqref{eq: k=n, alpha>0, bound on block P}, and the facts $\sqrt{\frac{d}{n}}<1$ and $\frac{\beta}{n^{\alpha}}<1$, for large $n$, we have

\begin{align}
    \expect{D} &\leq \frac{1}{d} + \frac{1}{n} \brac{1 + \frac{\beta}{n^\alpha} + O(\frac{1}{n^{2\alpha}})} \brac{1 + \sqrt{\frac{d}{n} + o(\frac{1}{n}}) + O(\frac{1}{n})}
    \nonumber
    \\ & \leq \frac{1}{d} + \frac{1}{n}+o(\frac{1}{n}).
\end{align}

The minimum achievable average response time of the accepted jobs is $\frac{1}{d}$. Therefore, the mean response time of the system converges to $\frac{1}{d}$ at rate $O(\frac{1}{n})$. 
\end{proof}



\section{Limited System Access} \label{section: Limited server access}

We now study the system with $k(n) < n$. Under this condition, the arrivals have access to a limited subset of servers of size $k(n)$, which are randomly sampled upon their arrival. We identify sufficient conditions on system parameters that guarantee the state space collapse results in the asymptotic limit. As a consequence, the system maintains the property of asymptotic optimality in terms of mean response time and blocking probability, even with limited subset access. However, the rate at which this optimality is achieved is different and depends on the specific characteristics of the sampled set.

In the following Lemma, we derive an upper bound on the probability of an arrival not finding $d$ idle servers in the sampled set of size $k(n)$.

\begin{lem} \label{lemma: zero assignment for k<n}
    If $q_1 \leq 1- \varepsilon$ for $\varepsilon >0$, we have
    $$1-A_d(x) \leq d (k(n))^d \brac{1-\varepsilon}^{k(n) - d}.$$
\end{lem}

\begin{proof}
    From the definition of assignment policies for $k(n) \leq n$ in Equation~\eqref{eq: arrival prob case k} we have

    \begin{equation*}
        1- A_d(x) = \sum_{i=0}^{d-1} \frac{\binom{nq_0}{i}\binom{nq_1}{k(n)-i}}{\binom{n}{k(n)}} =  \sum_{i=0}^{d-1} \binom{k(n)}{i} (q_0)^i (q_1)^{k(n)-i},
    \end{equation*}
    for $n$ sufficiently large. Moreover, from the fact $q_0 \leq 1$ and the assumption of Lemma $q_1 \leq 1- \varepsilon $ , we have

\begin{align*}
    1-A_d(x) &\leq \sum_{i=0}^{d-1} \frac{(k(n))^i (1-\varepsilon)^{k(n)-i}}{i!} \\&\leq d (k(n))^d (1-\varepsilon)^{k(n)-d}.
\end{align*}
for $n$ sufficiently large.
\end{proof}

\begin{lem} \label{lemma: first bound case k<n}
    If $\lambda(n)  = 1-\beta/n^{\alpha} \geq 0$ for $\alpha \in [0,1)$ and $\beta >0$, we have

    \begin{equation}
        \expect{(1-A_d(x))(\lambda (n) - dx_d)} \leq \lambda (n) d(k(n))^d \brac{1-\frac{\beta}{2n^{\alpha}}}^{k(n)-d} + d~ \expect{(\sum_{i \in [d]} x_i)\ind\brac{q_1 >1- \frac{\beta}{2n^{\alpha}}}}.
    \end{equation}
\end{lem}

\begin{proof}
We consider two cases where the fraction of busy servers exceeds the threshold $1- \frac{\beta}{2n^{\alpha}}$ and when it is below that threshold. We have

\begin{align}\label{eq: first term}
        \expect{(1-A_d(x))(\lambda (n)-dx_d)}
        =&\expect{(1-A_d(x))(\lambda (n)-dx_d)\ind\brac{q_1 \leq 1 - \frac{\beta}{2 n^{\alpha}}}} \\\label{eq: second term}+&
        \expect{(1-A_d(x))(\lambda (n)-dx_d)\ind\brac{q_1 > 1 - \frac{\beta}{2 n^{\alpha}}}} .
\end{align}

We bound each of the terms above in ~\eqref{eq: first term} and \eqref{eq: second term}.

Consider the term in expression~\eqref{eq: first term}. We assume $q_1 \leq 1- \frac{\beta}{2n^{\alpha}}$, otherwise this term will easily become zero due to the indicator function. Thus, we have $q_1 \leq 1- \varepsilon$ with $\varepsilon = \frac{\beta}{2n^{\alpha}}>0$, and from Lemma~\ref{lemma: zero assignment for k<n}, we conclude that 

\begin{equation}
    \expect{(1-A_d(x))(\lambda (n)-dx_d)\ind\brac{q_1 \leq 1- \frac{\beta}{2 n^\alpha}} } \leq \lambda(n) d (k(n))^d \brac{1- \frac{\beta}{2n^\alpha}}^ {k(n)-d}.
\end{equation} 

Consider the second term given by expression~\eqref{eq: second term}. We have

\begin{equation} \label{eq: main expression to bound}
    \begin{aligned}
        &\expect{(1-A_d(x))(\lambda (n)-dx_d)\ind\brac{q_1 > 1 - \frac{\beta}{2n^{\alpha}}}} \\ 
        \\ \leq& \expect{(1-A_d(x))(q_1- \frac{\beta}{2n^\alpha}-dx_d)\ind\brac{q_1 > 1 - \frac{\beta}{2n^{\alpha}}}}
        \\
        \leq&\expect{(1-A_d(x))(q_1-dx_d)\ind\brac{q_1 > 1 - \frac{\beta}{2n^{\alpha}}}}
        \\ =& \expect{(1-A_d(x))(\sum_{i \in [d-1]} ix_i)\ind\brac{q_1 > 1 - \frac{\beta}{2n^{\alpha}}}}
        \\ \leq & d~ \expect{(1-A_d(x))(\sum_{i \in [d-1]} x_i)\ind\brac{q_1 > 1 - \frac{\beta}{2n^{\alpha}}}}
        \\ \leq & d~\expect{(\sum_{i \in [d]} x_i)\ind\brac{q_1 > 1 - \frac{\beta}{2n^{\alpha}}}},
    \end{aligned}
\end{equation}

where the second line follows from the indicator function $\ind\brac{q_1 > 1- \frac{\beta}{2n^\alpha}}$ and the definition of  $\lambda(n) = 1- \frac{\beta}{n^\alpha}$; the fourth line follows from the definition of $q_1 = \sum_{i \in [d]}ix_i$; and the last line follows from the non-negativity of $x_i$. Combining all the bounds together, we get the desired result.
\end{proof}

In the following lemma, we establish an upper bound for expression~\eqref{eq: main expression to bound}. We introduce a new Lyapunov function $V_2(x) = (\sum_{i \in [d]} x_i) \ind \brac{q_1 > \lambda(n)+\delta}$ for a positive $\delta$, and show that outside of a suitable compact set, the drift of this Lyapunov function is strictly negative. Consequently, this implies that with high probability, the function $V_2(x)$ remains within that compact set.

Intuitively, when the fraction of busy servers exceeds a threshold that tends to one as $n$ grows, the number of accepted jobs into the system cannot increase substantially. In other words, the occurrence of two events of a significantly high number of busy servers and the acceptance of jobs into the system is highly improbable.

\begin{lem} \label{lemma: second Lyapunov function}
   For any $\delta \in (0 , \frac{\beta}{n^\alpha})$, define the following Lyapunov function.
    \begin{align}
    V_2(x) = (\sum_{i \in [d]} x_i)\ind\brac{q_1 > \lambda(n)+\delta}.
    \end{align}
If $V_2(x) \geq \kappa$ for some $\kappa>0$, then $GV_2(x) \leq - \delta $ and $\expect{V_2(x)} \leq \kappa + \frac{2}{n \delta}$, for all sufficiently large $n$.
    
\end{lem}
\begin{proof}
    Assume $V_2(x)  \geq \kappa$ for some $\kappa >0$. This implies
    \begin{equation}\label{eq: condition1}
        \ind\brac{q_1 > \lambda(n)+\delta} =1,
    \end{equation}
    and\begin{equation}\label{eq: condition2}
        q_1 > \lambda(n)+\delta.
    \end{equation}

    Let us calculate the drift of $V_2(x)$ under $G$, when $V_2(x)  \geq \kappa >0$. From the definition of the generator $G$ in Lemma~\ref{lemma: general policy lemma}, we have
\begin{equation*}
    GV_2(x) = \sum_{i \in [d]} n \lambda (n) A_i(x) \left(V_2(x+\frac{e_i}{n}) - V_2(x)\right)+ \sum_{i \in [d]} ni x_i \left(V_2(x-\frac{e_i}{n})-V_2(x)\right),
\end{equation*}
where $e_i$ denotes the d-dimensional unit vector with one at the $i^{th}$ position. Therefore,

\begin{equation*}
\begin{aligned}
    GV_2(x) &= \sum_{i \in [d]} n \lambda (n) A_i(x) \left(\brac{\sum_{j \in [d]} x_{j}+ \frac{e_i}{n}
    }\ind \brac{q_1+\frac{i}{n} > \lambda(n)+\delta} - (\sum_{j \in [d]} x_j) \ind\brac{q_1 > \lambda(n)+\delta}\right)
    \\&+ \sum_{i \in [d]} ni x_i \left(\brac{\sum_{j \in [d]} x_j - \frac{e_i}{n}}\ind\brac{q_1-\frac{i}{n} > \lambda(n)+\delta} - (\sum_{j \in [d]} x_j) \ind\brac{q_1 > \lambda(n)+\delta}\right).
\end{aligned}
\end{equation*}

However, under the assumption $V_2(x) \geq \kappa>0$, from Equations~\eqref{eq: condition1} and \eqref{eq: condition2}, we have

$$\ind\brac{q_1 + \frac{i}{n} > \lambda(n)+\delta} =1,\quad \text{for any } i \in[d],$$
and
$$\ind\brac{q_1 - \frac{i}{n} > \lambda(n)+\delta} = 1,\quad \text{for any } i \in[d], \text{ and } n \text{ sufficiently large.}$$

As a consequence, for $n$ large enough, we have
\begin{equation}
\begin{aligned}
    GV_2(x)&= \sum_{i \in [d]} n \lambda (n) A_i(x) (\frac{1}{n})+ \sum_{i \in [d]} ni x_i (-\frac{1}{n})
    \\
    &= \lambda (n) \sum_{i \in [d]} A_i(x)  - \sum_{i \in [d]} ix_i\\
    & \leq  \lambda (n) - q_1 \\
    & \leq -\delta,
\end{aligned}
\end{equation}

where the third line follows from the fact $\sum_{i \in [d]} A_i(x) \leq 1$ and the last line follows from \eqref{eq: condition2}. This shows that outside of the set $\{x: V_2(x) \leq \kappa, \kappa>0 \}$, the function has a negative drift and completes the first part of the lemma.

For the second part, we use the results from \cite[Theorem1 - (i)]{tail_bound_Gamarnik} which we recall in Appendix~\ref{appendix: tail_bound}. Under the notation of of Appendix~\ref{appendix: tail_bound}, we have $p_{max} = n$ and $\nu_{max}= 1/n$. Hence

\begin{align}
     \expect{V_2(x)} &\leq \kappa + \frac{2 p_{max}(\nu_{max})^2}{\delta}
     \\& \leq \kappa + \frac{2 n (1/n^2)}{\delta}
     \\& = \kappa + \frac{2}{n \delta}.
\end{align}
and this completes the proof.
\end{proof}

Combining Lemmas~\ref{lemma: general policy lemma},\ref{lemma: first bound case k<n} and~\ref{lemma: second Lyapunov function}, we obtain the following.

\begin{cor} \label{cor: k<n}
    If $\lambda(n)  = 1-\beta/n^{\alpha} \geq 0$ for $\alpha \in [0,1)$ and $\beta >0$, then we have
    
    \begin{equation}
    \begin{aligned}
        \expect{(\lambda (n) - dx_d)^2} \leq \frac{2d}{n}\lambda (n) \brac{1+\frac{2 n^{\alpha}}{\beta}}+ \lambda (n)^2 d(k(n))^d \brac{1-\frac{\beta}{2n^\alpha}}^{k(n)-d},
    \end{aligned}
    \end{equation}
    for all $n$ sufficiently large.
\end{cor}

\begin{proof}
    The result follows by choosing $\kappa = \frac{1}{n}$ and $\delta = \frac{\beta}{2n^\alpha}$ in Lemma~\ref{lemma: second Lyapunov function} and combining with Lemmas~\ref{lemma: general policy lemma} and \ref{lemma: first bound case k<n}.
\end{proof}

In the following theorem, we present the main performance bounds for systems with a finite size. We obtain sufficient conditions on the growth rate of the size of the sampling set $k(n)$, for the system to achieve asymptotic optimality. 

\begin{thm} \label{theorem: k<n, performance}
    Let $k(n)=n^\alpha \log(n)$ and $\lambda (n) = 1-\beta/n^{\alpha} \geq 0$, where $\alpha \in [0,1)$ and $\beta>0$. If $\beta >2(\alpha (d -1) +1)$ when $\alpha>0$, then in the steady state regime, the blocking probability of the system converges to zero at a rate of $O(n^{-\brac{1-\alpha}/2})$, and the mean response time of the accepted jobs converges to $\frac{1}{d}$ at a rate of $O(n^{-\brac{1-\alpha}/2})$.
\end{thm}    

\begin{proof}
    Let $P_b$ denote the blocking probability of the system. By the rate conservation law $\lambda(n)\brac{1-P_b} =\expect{q_1} \geq \expect{dx_d}$, and from the same arguments as in the proof of Theorem~\ref{theorem: k=n, performance}, we have 

    \begin{align} \label{eq: k<n, general bound on block}
        P_b^2 &\leq \frac{1}{\lambda(n)^2}\expect{\brac{\lambda(n) -dx_d}^2} 
       \nonumber\\
       & \leq  \frac{2d}{n\lambda (n)} \brac{1+\frac{2 n^{\alpha}}{\beta}}+ d n^{\alpha d } (\log(n))^d \brac{1-\frac{\beta}{2n^{\alpha}}}^{n^\alpha \log(n)-d},
    \end{align}
where the last inequality follows from Corollary~\ref{cor: k<n} when $k(n) = n^\alpha \log(n)$, and $n$ is sufficiently large. If $\alpha = 0$ and $\beta \in (0,1)$, we have

\begin{equation} \label{eq: block_P bound, k<n, alpha=0}
    P_b^2 \leq \frac{2d}{n(1-\beta)} \brac{1+\frac{2}{\beta}}+ o(\frac{1}{n}).
\end{equation}

Thus, the blocking probability of the system converges to zero at a rate of $O(\frac{1}{\sqrt{n}})$.

If $\alpha >0$ and $\beta > 2(\alpha d - \alpha +1)$, then from Equation~\eqref{eq: k<n, general bound on block} we have

\begin{align*}
    P_b^2 \leq \frac{2d}{n\lambda(n)} + \frac{4d}{\beta\lambda(n)}n^{\alpha-1}+ d n^{\alpha d} (\log(n))^d \brac{1-\frac{\beta}{2n^{\alpha}}}^{n^\alpha \log(n)-d}.
\end{align*}

Based on the properties of $\beta$, the second term in the above bound is the dominant term and we have

\begin{align*}
    P_b^2 &\leq \frac{4d}{\beta\lambda(n)}n^{\alpha-1} + o(n^{\alpha-1}).
\end{align*}

For $n$ large enough, $\frac{\beta}{n^\alpha} <1$, and we have

\begin{align} \label{eq: block_P bound, k<n, alpha>0}
    P_b^2 &\leq \frac{4d}{\beta}n^{\alpha-1} \brac{1+\frac{\beta}{n^\alpha} + O(\frac{1}{n^{2\alpha}})}+ o(n^{\alpha-1})
    \nonumber\\ &=\frac{4d}{\beta}n^{\alpha-1} + o(n^{\alpha-1}).
\end{align}

As a consequence, the blocking probability of the system converges to zero at a rate of $O(n^{-\brac{1-\alpha}/2})$.

Let $\expect{D}$ denote the mean response time of accepted jobs. From Little's law in the stationary regime, we have

\begin{equation*}
    \lambda(n) \brac{1-P_b} \expect{D} = \expect{\sum_{i \in d} x_i} = \expect{\sum_{i \in [d-1]} x_i} + \expect{x_d}.
\end{equation*}

For the first $d-1$ components of the system state, note that

$$\expect{\sum_{i \in [d-1]} x_i} \leq \expect{\sum_{i \in [d-1]} ix_i} = \expect{q_1 -dx_d} \leq \expect{\lambda(n) -dx_d},$$

where the equality follows from the definition of $q_1 = \sum_{i \in [d]} ix_i$; and the last inequality follows from the rate conservation law $\expect{q_1} = \lambda(n) (1-P_b) \leq \lambda(n)$. Moreover, for the $d^{th}$ component of the system state, we have $dx_d \leq q_1$ due to the definition of $q_1$ and $\expect{x_d} \leq \frac{\expect{q_1}}{d}$. Hence

\begin{align*}
    \expect{D} &\leq \frac{1}{\lambda(n) \brac{1- P_b}} \brac{\expect{\lambda(n) - dx_d} + \frac{\expect{q_1}}{d}}
     \\&= \frac{1}{\lambda(n) \brac{1- P_b}} \expect{\lambda(n) - dx_d} + \frac{1}{d},
\end{align*}

where the last line follows from the rate conservation law. Therefore,

\begin{align*}
    \brac{\expect{D} - \frac{1}{d}}^2 &\leq \frac{1}{\lambda(n)^2 \brac{1- P_b}^2} \brac{\expect{\lambda(n) - dx_d}}^2
    \\& \leq \frac{1}{\lambda(n)^2 \brac{1- P_b}^2} \expect{\brac{\lambda(n) - dx_d}^2},
\end{align*}

where the inequality follows from Jensen's inequality.

For $\alpha = 0$ and $\beta \in (0,1)$, from Corollary~\eqref{cor: k<n} and Equation~\eqref{eq: block_P bound, k<n, alpha=0}, we have

\begin{equation*}
    \expect{D - \frac{1}{d}}^2 \leq \frac{1}{\brac{1-\beta}^2 \brac{ 1-2\sqrt{\frac{2d}{n(1-\beta)}(1+\frac{2}{\beta}) + o(\frac{1}{n}})}} \brac{\frac{2d(1-\beta)}{n}\brac{1+\frac{2}{\beta}} + o(\frac{1}{n})}.
\end{equation*}

But $2\sqrt{\frac{2d}{n(1-\beta)}(1+\frac{2}{\beta}) + o(\frac{1}{n}}) <1$, for $n$ large enough and hence

\begin{align} \label{eq: response bound, k<n, alpha=0}
    \expect{D - \frac{1}{d}}^2 
    &\leq \frac{1}{(1-\beta)^2 } \brac{1+2\sqrt{\frac{2d}{n(1-\beta)}(1+\frac{2}{\beta}) + o(\frac{1}{n}})+ O(\frac{1}{n})}\brac{\frac{2d(1-\beta)}{n}\brac{1+\frac{2}{\beta}} + o(\frac{1}{n})}
    \nonumber
    \\& \leq \frac{2d}{n(1-\beta)}\brac{1+\frac{2}{\beta}}+ o(\frac{1}{n}).
\end{align}

This shows that the mean response time of accepted jobs converges to $\frac{1}{d}$ at a rate of $O(\frac{1}{\sqrt{n}})$.

If $\alpha>0$ and $\beta > 2(\alpha d - \alpha +1)$, from Corrolary~\ref{cor: k<n} and Equation~\eqref{eq: block_P bound, k<n, alpha>0} we have

\begin{equation*}
    \expect{D - \frac{1}{d}}^2 \leq \frac{1}{\brac{1- \frac{2\beta}{n^\alpha}} \brac{1-2\sqrt{\frac{4d}{\beta}n^{\alpha-1} + o(n^{\alpha-1}})}}\brac{\frac{4d}{\beta}n^{\alpha-1} + o(n^{\alpha-1})}.
\end{equation*}

For $n$ large enough, $\frac{2\beta}{n^\alpha}<1$ and $2\sqrt{\frac{4d}{\beta}n^{\alpha-1} + o(n^{\alpha-1}})<1$, and we have

\begin{align}
    \expect{D - \frac{1}{d}}^2 &\leq \brac{1+ \frac{2\beta}{n^\alpha} + O(\frac{1}{n^{2\alpha}})} \brac{1+ 2\sqrt{\frac{4d}{\beta}n^{\alpha-1} + o(n^{\alpha-1}})+ O(n^{\alpha-1})} \brac{\frac{4d}{\beta}n^{\alpha-1} +o(n^{\alpha-1})}
    \nonumber
    \\ &\leq \frac{4d}{\beta} n^{\alpha-1} + o(n^{\alpha-1}).
    \end{align}
As a consequence, the mean response time of accepted jobs converges to $\frac{1}{d}$ at a rate of $O(n^{-\brac{1-\alpha}/2})$, and this completes the proof.
\end{proof}

\section{Conclusion} \label{section: conclusion}
In this paper, we studied the performance of adaptive multiserver-job systems with $n$ servers. The system load was modeled as $1-\beta n^{-\alpha}\geq 0$, where $\alpha \in [0,1)$ and $\beta>0$. Upon arrival, each incoming job had the flexibility to split into a maximum of $d$ smaller tasks, based on the system state and available resources. We showed that when arriving jobs had complete knowledge of the system state and all servers were accessible, the system achieved the optimal average response time of $1/d$ and zero blocking probability in the limit as the system size approached infinity. We further characterized the rate of convergence, demonstrating a convergence rate of $O(1/\sqrt{n})$ for the blocking probability and $O(1/n)$ for the mean response time. When arrivals had partial knowledge of the system state, sampled upon their arrival, similar optimal performance measures were attainable. A necessary condition was that the sampling size must grow at rate $\omega\left(n^{\alpha}\right)$. In this case, the rate of convergence was established as $O(n^{-\brac{1-\alpha}/2})$ for both the mean delay and blocking probability.

The model we have introduced can be used to study systems where there exist job-server compatibility restrictions in the system and of interest in the context of content caching. In such systems, jobs can only be assigned to servers if they satisfy specific compatibility criteria. Incorporating these compatibility restrictions imposes challenges. Another possible avenue to explore is to consider threshold parallelization where the speedup before the threshold follows a concave function rather than a linear function as has been considered here. Addressing this question could yield valuable insights into system performance under more complex speedup dynamics.

\bibliographystyle{unsrt}
\bibliography{SteinRef.bib}

\begin{thebibliography}{10}

\bibitem{Fabrice_microservice}
Kiranpreet Kaur, Fabrice Guillemin, Veronica~Quintuna Rodriguez, and Francoise
  Sailhan.
\newblock Latency and network aware placement for cloud-native 5g/6g services.
\newblock In {\em 2022 IEEE 19th Annual Consumer Communications {\&} Networking
  Conference (CCNC)}, pages 114--119, Las Vegas, NV, USA, 2022. IEEE.

\bibitem{Google_Borg_2015}
Abhishek Verma, Luis Pedrosa, Madhukar Korupolu, David Oppenheimer, Eric Tune,
  and John Wilkes.
\newblock Large-scale cluster management at google with borg.
\newblock In {\em Proceedings of the Tenth European Conference on Computer
  Systems}, EuroSys '15, New York, NY, USA, 2015. Association for Computing
  Machinery.

\bibitem{MapReduce_2008}
Jeffrey Dean and Sanjay Ghemawat.
\newblock Mapreduce: Simplified data processing on large clusters.
\newblock {\em Commun. ACM}, 51(1):107–113, jan 2008.

\bibitem{TensorFlow_2016}
Mart{\'\i}n Abadi, Paul Barham, Jianmin Chen, Zhifeng Chen, Andy Davis, Jeffrey
  Dean, Matthieu Devin, Sanjay Ghemawat, Geoffrey Irving, Michael Isard,
  Manjunath Kudlur, Josh Levenberg, Rajat Monga, Sherry Moore, Derek~G. Murray,
  Benoit Steiner, Paul Tucker, Vijay Vasudevan, Pete Warden, Martin Wicke, Yuan
  Yu, and Xiaoqiang Zheng.
\newblock {TensorFlow}: A system for {Large-Scale} machine learning.
\newblock In {\em 12th USENIX Symposium on Operating Systems Design and
  Implementation (OSDI 16)}, pages 265--283, Savannah, GA, November 2016.
  USENIX Association.

\bibitem{ErasureCodes_2017}
Kangwook Lee, Nihar~B. Shah, Longbo Huang, and RamchandranKannan.
\newblock The mds queue: Analysing the latency performance of erasure codes.
\newblock {\em IEEE Transactions on Information Theory}, 63(5):2822--2842,
  2017.

\bibitem{moldable_2003}
Srinivasan, Krishnamoorthy, and Sadayappan.
\newblock A robust scheduling technology for moldable scheduling of parallel
  jobs.
\newblock In {\em 2003 Proceedings IEEE International Conference on Cluster
  Computing}, pages 92--99, Hong Kong, China, 2003. IEEE.

\bibitem{moldable_2002}
Walfredo Cirne and Francine Berman.
\newblock Using moldability to improve the performance of supercomputer jobs.
\newblock {\em Journal of Parallel and Distributed Computing},
  62(10):1571--1601, 2002.

\bibitem{Brill_multiserver_1984}
Percy~H. Brill and Linda Green.
\newblock Queues in which customers receive simultaneous service from a random
  number of servers: A system point approach.
\newblock {\em Management Science}, 30(1):51--68, 1984.

\bibitem{FILIPPOPOULOs_multiserver_2007}
Dimitrios Filippopoulos and Helen Karatza.
\newblock An m/m/2 parallel system model with pure space sharing among rigid
  jobs.
\newblock {\em Mathematical and Computer Modelling}, 45(5):491--530, 2007.

\bibitem{Zero_wait_central-2021}
Weina Wang, Qiaomin Xie, and Mor Harchol-Balter.
\newblock Zero queueing for multi-server jobs.
\newblock {\em SIGMETRICS Perform. Eval. Rev.}, 49(1):13–14, jun 2022.

\bibitem{sharp_zer_central_2022}
Yige Hong and Weina Wang.
\newblock Sharp waiting-time bounds for multiserver jobs.
\newblock In {\em Proceedings of the Twenty-Third International Symposium on
  Theory, Algorithmic Foundations, and Protocol Design for Mobile Networks and
  Mobile Computing}, MobiHoc '22, page 161–170, New York, NY, USA, 2022.
  Association for Computing Machinery.

\bibitem{Berg_speedup_2020}
Benjamin Berg, Mor Harchol-Balter, Benjamin Moseley, Weina Wang, and Justin
  Whitehouse.
\newblock Optimal resource allocation for elastic and inelastic jobs.
\newblock In {\em Proceedings of the 32nd ACM Symposium on Parallelism in
  Algorithms and Architectures}, SPAA '20, page 75–87, New York, NY, USA,
  2020. Association for Computing Machinery.

\bibitem{gamarnik2018delay}
David Gamarnik, John~N Tsitsiklis, and Martin Zubeldia.
\newblock Delay, memory, and messaging tradeoffs in distributed service
  systems.
\newblock {\em Stochastic Systems}, 8(1):45--74, 2018.

\bibitem{lu2011join}
Yi~Lu, Qiaomin Xie, Gabriel Kliot, Alan Geller, James~R Larus, and Albert
  Greenberg.
\newblock Join-idle-queue: A novel load balancing algorithm for dynamically
  scalable web services.
\newblock {\em Performance Evaluation}, 68(11):1056--1071, 2011.

\bibitem{Kaufman}
J.~Kaufman.
\newblock Blocking in a shared resource environment.
\newblock {\em IEEE Transactions on Communications}, 29(10):1474--1481, 1981.

\bibitem{Whitt_dropping_1985}
Ward Whitt.
\newblock Blocking when service is required from several facilities
  simultaneously.
\newblock {\em AT{\&}T Technical Journal}, 64(8):1807--1856, 1985.

\bibitem{liu2020steady}
Xin Liu and Lei Ying.
\newblock Steady-state analysis of load-balancing algorithms in the
  sub-halfin--whitt regime.
\newblock {\em Journal of Applied Probability}, 57(2):578--596, 2020.

\bibitem{braverman2020steady}
Anton Braverman.
\newblock Steady-state analysis of the join-the-shortest-queue model in the
  halfin--whitt regime.
\newblock {\em Mathematics of Operations Research}, 45(3):1069--1103, 2020.

\bibitem{FJ_two_servers_1984}
L.~Flatto and S.~Hahn.
\newblock Two parallel queues created by arrivals with two demands i.
\newblock {\em SIAM Journal on Applied Mathematics}, 44(5):1041--1053, 1984.

\bibitem{baccelli_FJ_1989}
François Baccelli, Armand~M. Makowski, and Adam Shwartz.
\newblock The fork-join queue and related systems with synchronization
  constraints: Stochastic ordering and computable bounds.
\newblock {\em Advances in Applied Probability}, 21(3):629–660, 1989.

\bibitem{Rizk_FJ_2016}
Amr Rizk, Felix Poloczek, and Florin Ciucu.
\newblock Stochastic bounds in fork–join queueing systems under full and
  partial mapping.
\newblock {\em Queueing Systems}, 83:261–291, 2016.

\bibitem{FJ_survey}
Alexander Thomasian.
\newblock Analysis of fork/join and related queueing systems.
\newblock {\em ACM Comput. Surv.}, 47(2), aug 2014.

\bibitem{heterogeneous_partial_FJ_2023}
Rooji Jinan, Gaurav Gautam, Parimal Parag, and Vaneet Aggarwal.
\newblock Asymptotic analysis of probabilistic scheduling for erasure-coded
  heterogeneous systems.
\newblock {\em SIGMETRICS Perform. Eval. Rev.}, 50(4):8–10, apr 2023.

\bibitem{limited_FJ_Wang_as_indep_2019}
Weina Wang, Mor Harchol-Balter, Haotian Jiang, Alan Scheller-Wolf, and
  R.~Srikant.
\newblock Delay asymptotics and bounds for multitask parallel jobs.
\newblock {\em Queueing Systems}, 46(3):207–239, jan 2019.

\bibitem{Wang_distributed_2020}
Wentao Weng and Weina Wang.
\newblock Achieving zero asymptotic queueing delay for parallel jobs.
\newblock {\em Proc. ACM Meas. Anal. Comput. Syst.}, 4(3):36, jun 2021.

\bibitem{multiserver_asymptotic_opimality_2022}
Isaac Grosof, Ziv Scully, Mor Harchol-Balter, and Alan Scheller-Wolf.
\newblock Optimal scheduling in the multiserver-job model under heavy traffic.
\newblock {\em Proc. ACM Meas. Anal. Comput. Syst.}, 6(3), dec 2022.

\bibitem{Amdahl's_law_2008}
Mark~D. Hill and Michael~R. Marty.
\newblock Amdahl's law in the multicore era.
\newblock {\em Computer}, 41(7):33--38, 2008.

\bibitem{Berg_speedup_2018}
Benjamin Berg, Jan-Pieter Dorsman, and Mor Harchol-Balter.
\newblock Towards optimality in parallel job scheduling.
\newblock {\em SIGMETRICS Perform. Eval. Rev.}, 46(1):116–118, jun 2018.

\bibitem{WCFS_Grososf_2022}
Isaac Grosof, Mor Harchol-Balter, and Alan Scheller-Wolf.
\newblock Wcfs: A new framework for analyzing multiserver systems.
\newblock {\em Queueing Systems}, 102(1-2):143--174, 2022.

\bibitem{tail_bound_Gamarnik}
Dimitris Bertsimas, David Gamarnik, and John~N Tsitsiklis.
\newblock Performance of multiclass markovian queueing networks via piecewise
  linear lyapunov functions.
\newblock {\em Annals of Applied Probability}, 11(4):1384--1428, 2001.

\end{thebibliography}
\appendix

\section{Tail Bounds}
\label{appendix: tail_bound}
We recall in this section Theorem $1$ of \cite{tail_bound_Gamarnik}. Consider a continuous time Markov chain $X(t)$ which takes values in some countable set $\mathcal{X}$, with a stationary probability distribution $\pi$. For any two vectors $x, x^' \in \mathcal{X}$, let $p(x, x^')$ denote the transition probability from state $x$ to state $x^'$. For a given function $\Phi : \mathcal{X} \to \mathbb{R}^+$ such that $\mathbb{E}_{\pi}\left[\Phi(X(t))\right]< \infty$, let 
$$\nu_{max}= \underset{x,x^' \in \mathcal{X}, p(x,x^')>0} {\sup} |\Phi(x')-\Phi(x)|< \infty$$
Namely, $\nu_{max}$ is the largest possible change of the function $\Phi$ during an arbitrary transition. Also, let 
$$p_{max} = \underset{x \in \mathcal{X}}{\sup} \sum_{x^' \in \mathcal{X}, \Phi(x) < \Phi(x^')}p(x,x^')<\infty.$$
Namely, $p_{max}$is the tight upper bound on the probability that the value of $\Phi$ is increasing during an arbitrary transition.

\begin{enumerate}[(i)]
    \item If there exists a Lyapunov function $\Phi$ such that for any $x \in \mathcal{X}$ with $\Phi(x) > B$, 
    $$G \Phi (x) = \sum_{x^' \not= x} p(x,x^') \brac{\Phi(x^') - \Phi(x)} \leq - \gamma,$$
for some $\gamma > 0$ and $B \geq 0$, then for any $m = 0, 1, 2, ...,$
    
    \begin{equation*}
        \mathbb{P}_{\pi}\left(\Phi(X(t))> B + 2\nu_{max}m \right) \leq \left(\frac{p_{max} \nu_{max}}{p_{max} \nu_{max} + \gamma}\right)^{m+1}.
    \end{equation*}
    As a result,

    \begin{equation*}
       \mathbb{E}_{\pi}\left[\Phi(X(t))\right] \leq B + \frac{2p_{max} (\nu_{max})^2}{\gamma}.
    \end{equation*}
\end{enumerate}
\end{document}